\documentclass[12pt]{article}

\usepackage[numbers]{natbib}
\usepackage{rotating}
\usepackage{graphicx}

\begin{document}

\title{Can We Rely on AI?}

\author{
Desmond J. Higham\thanks{School of Mathematics, University of Edinburgh, James Clerk Maxwell Building, Edinburgh, UK. This manuscript was prepared as an Extended Abstract submission for the 21st International Conference of Numerical Analysis and Applied Mathematics,
Heraklion, Crete, Greece, 2023.
} }

\date{August 2023}

\maketitle

\begin{abstract}
Over the last decade, adversarial attack algorithms have  
revealed instabilities in deep learning tools. These algorithms raise issues 
regarding safety, reliability and interpretability in artificial intelligence; especially in high risk settings. From a practical perspective, there has been 
a war of escalation between those developing attack and defence strategies. 
 At a more theoretical level, researchers have also studied bigger picture questions concerning the existence and computability of attacks. 
 Here we give a brief overview of the topic, focusing on aspects that are likely to be 
 of interest to researchers in applied and computational mathematics.
\end{abstract}

\section{Introduction}\label{sec:intro}
We are currently living in a world where 
\begin{itemize}
\item 
automated driving systems can misinterpret
 traffic 
 ``Stop'' signs as speed limit signs when  minimal graffiti is added
\cite{physical},
\item artificial intelligence (AI) algorithms
for medical imaging, microscopy and other inverse problems in the
sciences can introduce hallucinations \cite{Hansen_PNAS_2020,Hallucination_NatureM2},
\item methods that aim to ``explain'' or ``interpret'' the decisions from an AI system
may be superficial or unreliable \cite{NeurIPSSaliency,falsehope21}, 
\item carefully tailored prompts can persuade 
chatbox-style large language models to generate offensive or sensitive material 
\cite{BLJCHD2323,CNCJGAWILTS23,WWALS23}.
\end{itemize}
 
 At the heart of these issues is the concept of \emph{instability}---extreme sensitivity of the output to changes in the input. In this brief manuscript we 
 describe how AI instabilities have been identified and analysed,
aiming at a readership in  
 applied and computational mathematics.

 \section{Adversarial Attacks}\label{sec:adv}
 
\begin{figure}[h]
  \includegraphics[height=180pt]{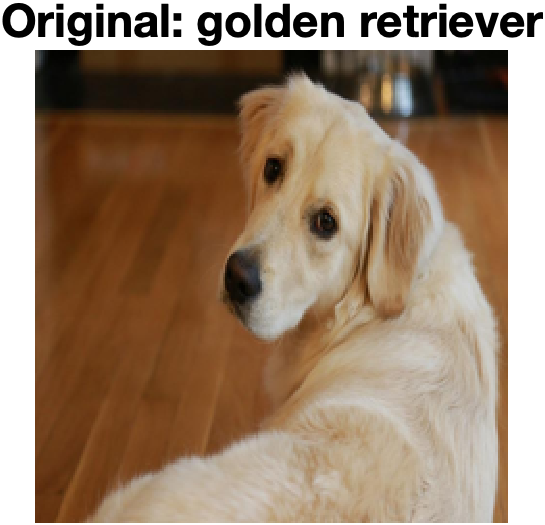} 
  \quad
  \includegraphics[height=180pt]{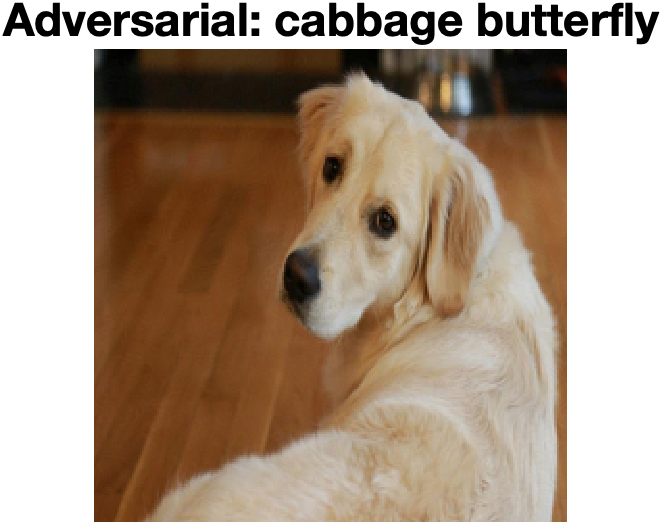}
  
  \caption{Left: an image from ImageNet that is correctly classified as a golden retriever by 
  a convolutional neural network. Right: after a small,
  visually imperceptible, perturbation is made, the image is now classified as a cabbage butterfly. 
  \label{fig:adv}
 }
\end{figure}

 Let us focus now on image classification. 
 In 
 Figure~\ref{fig:adv}
  we show a successful  \emph{adversarial attack}.
  This computation was performed with the MATLAB Deep Learning Toolbox
  \cite{MATLAB:2022},
  after slightly editing the relevant demonstration code\footnote{
  \texttt{https://uk.mathworks.com/help/deeplearning/ug/generate-adversarial-examples.html}.
  }.
  (The only significant edit was changing ``great white shark'' to ``cabbage butterfly.'')
  The computation uses \texttt{squeezenet}, an 18 layer
  convolutional neural network that has been trained on over a million 
  images from the ImageNet database \texttt{http://www.image-net.org}.
  The network classifies an image into  one of 1000 object categories.
  On the left in 
 Figure~\ref{fig:adv}
 we show an image from the ImageNet data set which is correctly classified 
 by \texttt{squeezenet} as a golden retriever.
 On the right we show the image arising when a specific, small perturbation is made.
  The perturbation was computed using an attack algorithm that has access to 
  the inner workings of the network, along the lines of the ideas in 
   \cite{kurakin2016adversarial,towards19}.  Here, the 
   attack algorithm is designed to produce a perturbed image that is unchanged
    according to the human eye, but is now classified 
   by \texttt{squeezenet}     
    as a cabbage butterfly, and 
    these goals were achieved.
    
    The observation that convolutional neural networks are vulnerable 
    to computable adversarial perturbations was made in 
    \cite{szegedy2013intriguing}; see also 
    \cite{harness} for further seminal work on this topic.
    Such perturbations may be found by evaluating, or estimating, the 
    sensitivity of the classifier output (or the loss function that is used to 
    judge goodness-of-fit) to small changes in the input pixel values. A ``steepest ascent'' direction is then available that, locally,
    has the most damaging effect. 
    
    Many variants of this basic idea have been studied:
    \begin{itemize}
     \item In Figure~\ref{fig:adv} we formed a \emph{targeted attack}, asking for the new class to be ``cabbage butterfly.'' An \emph{untargeted attack} would have the simpler aim of 
     seeking any change to the output class.
     \item Figure~\ref{fig:adv} shows a \emph{white box} attack---the attack algorithm has access to the partial derivatives of the classifier output and loss function with respect to the input pixels. A corresponding \emph{black box} attack has access only to input and output pairs; in this case certain partial derivatives can be estimated via finite differencing.
     \item The attack in Figure~\ref{fig:adv} applies to a single image.
     In practice it has been observed that the same perturbation can lead to misclassification 
     across many images and across different deep learning networks, leading to so-called
     \emph{universal attacks} \cite{moosavi17,ko-uapsvd-2018}.
      \item The algorithm in Figure~\ref{fig:adv} looked for a small perturbation as measured  in the Euclidean norm (regarding the pixels as forming a large input vector).
      Alternative norms may be used 
       \cite{HXSS15,onep}
        and more tailored relative componentwise perturbations may be implemented
        \cite{BH23}.
       \item Attacks may be made on the parameters that define the network
       \cite{breier18,badnets19,hammer20,subnet21,tyukin2020adversarial,thwg21}.
           \end{itemize} 
     
     In tandem with attack algorithms, there have been a wealth of defence strategies; 
      including \emph{adversarial training} 
      \cite{Madry18}, where 
      a network is retrained on perturbed data points using a robust optimization approach.
     Currently, within this extensive      
     cat-and-mouse game there is no guaranteed safeguard against adversaries  \cite{Attack_survey_2018}.
       Indeed, the author in \cite{Car23} argues that ``Historically, the vast majority of adversarial defenses published at top-tier conferences \ldots are quickly broken.''
       We also note that vulnerability to misclassification 
       raises major questions around the usefulness of add-on strategies that aim to ``explain''
       or ``interpret'' system outputs
       \cite{NeurIPSSaliency,falsehope21}.
       For the picture on the right in Figure~\ref{fig:adv}, there is clearly no valid 
       post-hoc explanation for the ``cabbage butterfly'' classification.

 \section{Higher Level Issues}\label{sec:inev}
 
Alongside the development of heurustic attack and defence algorithms, researchers have also
addressed the overarching question of whether, under reasonable assumptions, it is inevitable that classification tools will be vulnerable to attack.
Of course, in order to prove rigorous results it is necessary to make assumptions about the 
training data, testing data and classifier.
In \cite{BHV21} the accuracy-robustness tradeoff is studied, and it is shown that any approach to 
training a neural network with a fixed architecture must allow for examples of either inaccuracy or instability.
The authors in \cite{shafahi2018adv,tyukin2020adversarial} look at the asymptotic limit of 
high-dimensional data (i.e., large numbers of pixels),
where concentration of measure effects are relevant, and show that under appropriate assumptions it is inevitable that adversarial examples will exist.
In \cite{thwg21} a concrete ``stealth attack'' algorithm 
that perturbs parameters in a neural network
is introduced. Unusually, this algorithm is amenable to rigorous analysis, and it can be shown to have a high probability of success in an appropriate high-dimensional limit.

Although ``inevitability of vulnerability'' results paint a worrying picture for the 
use of AI, they may also help in the design of more stable systems. In order to 
work around an inevitability result we must ensure that one or more of the 
assumptions under which the result was proved does not hold. 
In the case of results that apply in asymptotically high dimension, one reasonable approach is to estimate, and, if necessary reduce, the real, or \emph{intrinsic}, dimension
that operates in the input space and through the layers of a neural network---this might be 
much less than the nominal dimension \cite{ZG23}.
(For example, the space of natural images may lie on a lower dimensional manifold than the overall pixel space.) 
We also note that concentration of measure effects can be used to motivate on-the-fly fixes 
for large, expensively trained systems when a small number of errors have been identified
\cite{Tyukin_Notices_2022}.
Moreover, there has been some progress in establishing that 
guaranteed benefits arise from  
adversarial training; see, for example, 
 \cite{BS22}.

 \section{Implications for Regulation}\label{sec:reg}  
 
 The increasing prevalence of AI in our everyday lives is raising concerns around
 data protection, 
 privacy, safety, ethics, defence, security, and the use of the earth's natural resources.
 Understanding and addressing these issues requires a multi-disciplinary approach,
 combining expertise across the computational, engineering, political and social sciences, 
  as well as 
  business, law and philosophy.
  We argue here that more fundamental concerns, arising   
    at a mathematical level, should be added to this mix.
  AI systems may possess inherent vulnerabilities, preventing them from 
  operating as intended and making them susceptible to adversarial attacks. 
  
  These concerns are becoming highly relevant as governmental bodies seek to 
  regulate AI.   
  For example:  
 ``The AI Act is a proposed regulation by the European Union that aims to establish a legal framework for the development, deployment, and use of AI systems in the EU'' \cite{eu23}.
The current version of the proposed act would oblige providers of AI systems 
intended for high-risk applications to comply with a range of requirements,
including \emph{technical robustness}.
Similarly,  
 Section~3.2.3 of 
 the 
 UK Government policy paper \emph{AI regulation: a pro-innovation approach} (updated August, 2023)
 recommends five cross-sectoral principles, two of which are
  \emph{Safety, security and robustness} and 
 \emph{Appropriate transparency and explainability}.
 There is a well-defined sense in which such 
 goals are 
 mathematically unachievable.
 In order to create regulations that can withstand mathematical scrutiny, 
 it is clear that
  terms such as robustness, safety and explainability 
 must be articulated in sufficient detail (along with the term AI itself), and known results about the 
 fundamental limitations 
 of algorithms that are designed to perform inference and 
 decision making must be acknowledged. 
 
 \bigskip

\noindent
\textbf{Acknowledgements}
DJH was supported 
  by EPSRC grants EP/P020720/1 and EP/V046527/1. 
  The experiment summarized in Figure~\ref{fig:adv} can be performed 
  by changing ``great white shark'' to ``cabbage butterfly'' in  
the demonstration code at\\
{\scriptsize
  \texttt{https://uk.mathworks.com/help/deeplearning/ug/generate-adversarial-examples.html}.
  }



\bibliographystyle{siam}
\bibliography{adv_refs}

\begin{thebibliography}{10}

\bibitem{towards19}
{\em Towards deep learning models resistant to adversarial attacks}, arXiv
  preprint arXiv: 1706.06083,  (2019).

\bibitem{Attack_survey_2018}
{\sc N.~Akhtar and A.~Mian}, {\em Threat of adversarial attacks on deep
  learning in computer vision: A survey}, IEEE Access, 6 (2018),
  pp.~14410--14430.

\bibitem{Hansen_PNAS_2020}
{\sc V.~Antun, F.~Renna, C.~Poon, B.~Adcock, and A.~C. Hansen}, {\em On
  instabilities of deep learning in image reconstruction and the potential
  costs of {AI}}, Proc. Natl. Acad. Sci. USA,  (2020).

\bibitem{BHV21}
{\sc A.~Bastounis, A.~C. Hansen, and V.~Vla\^ci\'c}, {\em The mathematics of
  adversarial attacks in {AI}--{W}hy deep learning is unstable despite the
  existence of stable neural networks}, arXiv:2109.06098 [cs.LG],  (2021).

\bibitem{BH23}
{\sc L.~Beerens and D.~J. Higham}, {\em {Adversarial ink: componentwise
  backward error attacks on deep learning}}, IMA Journal of Applied
  Mathematics,  (2023), p.~hxad017.

\bibitem{breier18}
{\sc J.~Breier, X.~Hou, D.~Jap, L.~Ma, S.~Bhasin, and Y.~Liu}, {\em Practical
  fault attack on deep neural networks}, in ACM SIGSAC Conference on Computer
  and Communications Security (CCS), ACM, 2018, pp.~2204--2206.

\bibitem{BS22}
{\sc L.~Bungert and K.~Stinson}, {\em Gamma-convergence of a nonlocal perimeter
  arising in adversarial machine learning}, arXiv preprint arXiv: 2211.15223,
  (2022).

\bibitem{Car23}
{\sc N.~Carlini}, {\em A {LLM} assisted exploitation of {AI}-{G}uardian}, arXiv
  preprint arXiv:2307.15008,  (2023).

\bibitem{CNCJGAWILTS23}
{\sc N.~Carlini, M.~Nasr, C.~A. Choquette-Choo, M.~Jagielski, I.~Gao,
  A.~Awadalla, P.~W. Koh, D.~Ippolito, K.~Lee, F.~Tram{\`e}r, and L.~Schmidt},
  {\em Are aligned neural networks adversarially aligned?}, arXiv preprint
  arXiv:2306.15447,  (2023).

\bibitem{ZG23}
{\sc Z.~Cui and P.~Grindrod}, {\em {Mappings, dimensionality and reversing out
  of deep neural networks}}, IMA Journal of Applied Mathematics,  (2023),
  p.~hxad019.

\bibitem{NeurIPSSaliency}
{\sc A.-K. Dombrowski, M.~Alber, C.~Anders, M.~Ackermann, K.-R. M\"uller, and
  P.~Kessel}, {\em Explanations can be manipulated and geometry is to blame},
  in NeurIPS, H.~Wallach, H.~Larochelle, A.~Beygelzimer, F.~d'~Alch\'e-Buc,
  E.~Fox, and R.~Garnett, eds., vol.~32, Curran Associates, Inc., 2019.

\bibitem{physical}
{\sc I.~Evtimov, K.~Eykholt, E.~Fernandes, T.~Kohno, B.~Li, A.~Prakash,
  A.~Rahmati, and D.~Song}, {\em Robust physical-world attacks on machine
  learning models}, CoRR, abs/1707.08945 (2017).

\bibitem{falsehope21}
{\sc M.~Ghassemi, L.~Oakden-Rayner, and A.~L. Beam}, {\em The false hope of
  current approaches to explainable artificial intelligence in health care},
  The Lancet Digital Health, 3 (2021), pp.~e745--e750.

\bibitem{harness}
{\sc I.~J. Goodfellow, J.~Shlens, and C.~Szegedy}, {\em Explaining and
  harnessing adversarial examples}, in 3rd International Conference on Learning
  Representations, San Diego, CA, Y.~Bengio and Y.~LeCun, eds., 2015.

\bibitem{Tyukin_Notices_2022}
{\sc A.~Gorban, B.~Grechuk, and I.~Tyukin}, {\em Stochastic separation
  theorems: How geometry may help to correct {AI} errors}, Notices Amer. Math.
  Soc., 70 (2022), pp.~25--33.

\bibitem{badnets19}
{\sc T.~Gu, K.~Liu, B.~Dolan-Gavitt, and S.~Garg}, {\em {B}ad{N}ets:
  {E}valuating backdooring attacks on deep neural networks}, IEEE Access, 7
  (2019), pp.~47230--47244.

\bibitem{Hallucination_NatureM2}
{\sc D.~P. Hoffman, I.~Slavitt, and C.~A. Fitzpatrick}, {\em The promise and
  peril of deep learning in microscopy}, Nature Methods, 18 (2021),
  pp.~131--132.

\bibitem{HXSS15}
{\sc R.~Huang, B.~Xu, D.~Schuurmans, and C.~Szepesvari}, {\em Learning with a
  strong adversary}, arXiv preprint arXiv:1511.03034,  (2017).

\bibitem{ko-uapsvd-2018}
{\sc V.~Khrulkov and I.~Oseledets}, {\em Art of singular vectors and universal
  adversarial perturbations}, in IEEE Conference on Computer Vision and Pattern
  Recognition (CVPR), 2018.

\bibitem{kurakin2016adversarial}
{\sc A.~Kurakin, I.~Goodfellow, and S.~Bengio}, {\em Adversarial examples in
  the physical world}, arXiv preprint arXiv:1607.02533,  (2016).

\bibitem{BLJCHD2323}
{\sc B.~Liu, B.~Xiao, X.~Jiang, S.~Cen, X.~He, and W.~Dou}, {\em Adversarial
  attacks on large language model-based system and mitigating strategies: {A}
  case study on {C}hat{GPT}}, Security and Communication Networks,  (2023),
  p.~8691095.

\bibitem{eu23}
{\sc T.~Madiega}, {\em Briefing {D}ocument: {A}rtificial intelligence act,
  {E}uropean {P}arliamentary {R}esearch {S}ervice}, 2023.

\bibitem{Madry18}
{\sc A.~Madry, A.~Makelov, L.~Schmidt, D.~Tsipras, and A.~Vladu}, {\em Towards
  deep learning models resistant to adversarial attacks}, in 6th International
  Conference on Learning Representations, Vancounver, 2018.

\bibitem{MATLAB:2022}
{\sc MATLAB}, {\em version 9.13.0.2080170 (R2022b)}, The MathWorks Inc.,
  Natick, Massachusetts, 2022.

\bibitem{moosavi17}
{\sc S.~Moosavi-Dezfooli, A.~Fawzi, O.~Fawzi, and P.~Frossard}, {\em Universal
  adversarial perturbations}, in 2017 IEEE Conference on Computer Vision and
  Pattern Recognition (CVPR), 2017, pp.~86--94.

\bibitem{subnet21}
{\sc X.~Qi, J.~Zhu, C.~Xie, and Y.~Yang}, {\em Subnet replacement:
  {D}eployment-stage backdoor attack against deep neural networks in gray-box
  setting}, ICLR 2021 Workshop on Security and Safety in Machine Learning
  System,  (2021).

\bibitem{shafahi2018adv}
{\sc A.~Shafahi, W.~Huang, C.~Studer, S.~Feizi, and T.~Goldstein}, {\em Are
  adversarial examples inevitable?}, International Conference on Learning
  Representations, New Orleans, USA,  (2019).

\bibitem{onep}
{\sc J.~Su, D.~V. Vargas, and S.~Kouichi}, {\em One pixel attack for fooling
  deep neural networks}, arXiv preprint arXiv: 1710.08864,  (2017).

\bibitem{szegedy2013intriguing}
{\sc C.~Szegedy, W.~Zaremba, I.~Sutskever, J.~Bruna, D.~Erhan, I.~Goodfellow,
  and R.~Fergus}, {\em Intriguing properties of neural networks}, arXiv
  preprint arXiv:1312.6199,  (2013).

\bibitem{thwg21}
{\sc I.~Y. Tyukin, D.~J. Higham, A.~Bastounis, E.~Woldegeorgis, and A.~N.
  Gorban}, {\em The feasibility and inevitability of stealth attacks},
  arXiv:2106.13997,  (2021).

\bibitem{tyukin2020adversarial}
{\sc I.~Y. Tyukin, D.~J. Higham, and A.~N. Gorban}, {\em On adversarial
  examples and stealth attacks in artificial intelligence systems}, in 2020
  International Joint Conference on Neural Networks, IEEE, 2020, pp.~1--6.

\bibitem{WWALS23}
{\sc Y.~Wolf, N.~Wies, O.~Avnery, Y.~Levine, and A.~Shashua}, {\em Fundamental
  limitations of alignment in large language models}, arXiv preprint
  arXiv:2304.11082,  (2023).

\bibitem{hammer20}
{\sc F.~Yao, A.~S. Rakin, and D.~Fan}, {\em Deep{H}ammer: {D}epleting the
  intelligence of deep neural networks through targeted chain of bit flips},
  Proceedings of the 29th USENIX Security Symposium,  (2020).

\end{thebibliography}

\end{document}